\documentclass[12pt]{amsart}
\usepackage{amsmath, amsthm, graphicx, amssymb,verbatim,pst-all,color,txfonts}
\usepackage[all]{xy}
\usepackage[margin=1.0in]{geometry}
\ifx\JPicScale\undefined\def\JPicScale{1.0}\fi
\unitlength \JPicScale mm
{\theoremstyle{plain}
   
   \newtheorem{theorem}{Theorem}[section]

}
{\theoremstyle{definition}

   \newtheorem{remark}[theorem]{Remark}

}

\newcommand{\M}{\overline{\mathcal{M}}}

\newcommand{\rk}{\operatorname{rk}}

\numberwithin{theorem}{section}

\begin{document}
\title{Fibonacci, golden ratio, and vector bundles}
\author{Noah Giansiracusa}
\maketitle

%\tableofcontents

\begin{abstract}
There is a family of vector bundles over the moduli space of stable curves that, while first appearing in theoretical physics, has been an active topic of study for algebraic geometers since the 1990s.  By computing the rank of the exceptional group $G_2$ case of these bundles in three different ways, we derive a family of summation formulas for Fibonacci numbers in terms of the golden ratio.
\end{abstract}

\section{Introduction}

The purpose of this note is to establish the following formulas:
\begin{theorem}\label{thm:main}
For all integers $n \ge 0$ and $g \ge 2$, 
\begin{equation}\label{eq:main}
\sum_{i=0}^g \binom{g}{i}F_{2i+n-1} = \varphi^n(\varphi+2)^{g-1} + \overline{\varphi}^n(\overline{\varphi}+2)^{g-1} = \sum_{i=0}^g \binom{g}{i}2^{g-i}F_{i+n-1},
\end{equation}
where $F_k$ is the $k^{th}$ Fibonacci number, $\varphi = \frac{1+\sqrt{5}}{2}$ is the golden ratio and $\overline{\varphi} = \frac{1-\sqrt{5}}{2}$ is the conjugate root of its minimal polynomial $x^2 - x -1$.  
\end{theorem}
While these formulas could be proven using elementary means, our emphasis in this paper is how naturally they emerge from an ostensibly unrelated setting.  Indeed, by computing the rank of a certain vector bundle, for each $n$ and $g$, in three different ways we arrive at the above formulas.  

The relevant vector bundles here are the exceptional Lie algebra $\mathfrak{g}_2$ case of a more general construction, occurring for each simple complex Lie algebra, of the so-called vector bundle of \emph{conformal blocks} on the Deligne-Mumford moduli space $\M_{g,n}$ of stable $n$-pointed genus $g$ curves.  These vector bundles first appeared in mathematical physics \cite{KZ84,Ver88} and they have been studied from the perspectives of both infinite-dimensional representation theory \cite{TUY89,Uen97,Loo05} and algebraic geometry \cite{NR93,BL94,Fal94,Pau96,Fak12,MOP15}; these citations are just a sample of the vast literature.  One of the most significant properties of these bundles is their recursive structure reflecting the boundary stratification of $\M_{g,n}$.  This property, known as \emph{factorization}, was cleverly used by Mukhopadhyay to compute the rank of the $\mathfrak{g}_2$, level 1, case of these bundles and the formula he obtained \cite{Muk16} is the middle expression in Equation \eqref{eq:main} involving the golden ratio; by using factorization in two different ways we obtain the two Fibonacci summation expressions in that equation.  The details of this assertion are the content of this paper.  

\begin{remark}
The author learned from the On-Line Encyclopedia of Integer Sequences that in the case $n=0$, the formulas in Theorem \ref{thm:main} count the number of no-leaf edge-subgraphs in the M\"obius ladder $M_{g-1}$ \cite{OEIS} (see \cite[Theorem 16.1]{McS98}).\end{remark}

\subsection*{Acknowledgements}
The author was supported in part by NSF DMS-1802263 and thanks Swarnava Mukhopadhyay for generously sharing ideas that led to a substantial revision of an earlier version of this paper.

\section{Conformal blocks vector bundles}

\subsection{Notation}

Given a simple complex Lie algebra $\mathfrak{g}$, an integer $\ell \ge 1$ called the \emph{level}, and an $n$-tuple \[\underline{\mu}=(\mu_1,\ldots,\mu_n) \in P_\ell\] of dominant integral weights for $\mathfrak{g}$ of level $\ell$, there is a vector bundle of \emph{conformal blocks} on the moduli stack $\M_{g,n}$ of stable $n$-pointed genus $g$ curves that we denote by $\mathbb{V}_{g,n}^{\mathfrak{g},\ell}(\underline{\mu})$, or by $\mathbb{V}^{\mathfrak{g},\ell}_g$ when $n=0$ \cite{TUY89,Uen97,Fak12}.  The construction of this vector bundle is immaterial for our purposes, we shall only need some standard properties well-established in the literature.  

\subsection{The Verlinde formula}

The main property of conformal blocks vector bundles we are concerned with is their rank.  There is a famous direct formula for this, the \emph{Verlinde formula} (see \cite{Sor96} for a nice survey), which for $n=0$ and $\mathfrak{g}$ of classical type or $\mathfrak{g}_2$ reads as follows (see \cite[Proposition B.1]{Gre11}):
\begin{equation}\label{eq:Verlinde}
\rk \mathbb{V}^{\mathfrak{g},\ell}_g = \left((\ell+h^\vee)^{\rk(\mathfrak{g})}\#(P/Q)\#(Q/Q_{lg})\right)^{g-1}\sum_{\mu\in P_\ell}\prod_{\alpha\in\Delta_+}2\sin\left(\frac{\pi\langle \alpha, \mu+\rho\rangle}{\ell+h^\vee}\right),
\end{equation}
where $P_\ell$ is the set of dominant integral weights for $\mathfrak{g}$ of level $\ell$, $\Delta_+$ is the set of positive roots, $\langle \cdot,\cdot \rangle$ is the Killing form, $\rho = \frac{1}{2}\sum_{\alpha_j\in\Delta_+}\alpha_j$, $h^\vee$ is the dual Coxeter number, $P$ is the weight lattice, $Q$ is the root lattice, and $Q_{lg}$ is the long root lattice.  

Gregoire has evaluated this explicitly in the case that is relevant for this paper, $\mathfrak{g}=\mathfrak{g}_2$ and $\ell=1$ \cite[Appendix B]{Gre11} (see also \cite{GP13}):
\begin{equation}\label{eq:Gregoire}
\rk \mathbb{V}^{\mathfrak{g}_2,1}_g = \left(\frac{5+\sqrt{5}}{2}\right)^{g-1} + ~ \left(\frac{5-\sqrt{5}}{2}\right)^{g-1}
\end{equation}
This is the middle expression in Equation \eqref{eq:main} when $n=0$.  

\subsection{Factorization, fusion rules, propagation of vacua}

Another way to compute the rank of conformal blocks vector bundles is via \emph{factorization}, which describes the behavior of the fibers over a curve as the curve degenerates from smooth to nodal \cite{TUY89}.  Fakhruddin has nicely expressed this in terms of boundary divisors in the moduli space \cite[Proposition 2.4]{Fak12}: if
\[\gamma : \M_{g_1,n_1+1} \times \M_{g_2,n_2+1} \rightarrow \M_{g_1+g_2,n_1+n_2}\] is the gluing map attaching the pair of points given by the last marked point on each component, then 
\begin{equation*}
\gamma^*\mathbb{V}^{\mathfrak{g},\ell}_{g_1+g_2,n_1+n_2}(\mu_1,\ldots,\mu_{n_1+n_2}) \cong \bigoplus_{\mu\in P_\ell} \mathbb{V}^{\mathfrak{g},\ell}_{g_1,n_1+1}(\mu_1,\ldots,\mu_{n_1},\mu) \otimes \mathbb{V}^{\mathfrak{g},\ell}_{g_2,n_2+1}(\mu_{n_1+1},\ldots,\mu_{n_1+n_2},\mu^*).
\end{equation*}
Taking ranks, this says
\begin{equation}\label{eq:fact1}
\rk\mathbb{V}^{\mathfrak{g},\ell}_{g_1+g_2,n_1+n_2}(\mu_1,\ldots,\mu_{n_1+n_2}) = \sum_{\mu\in P_\ell} \rk \mathbb{V}^{\mathfrak{g},\ell}_{g_1,n_1+1}(\mu_1,\ldots,\mu_{n_1},\mu) \rk \mathbb{V}^{\mathfrak{g},\ell}_{g_2,n_2+1}(\mu_{n_1+1},\ldots,\mu_{n_1+n_2},\mu^*)
\end{equation}
Similarly, if
\[\gamma : \M_{g-1,n+2} \rightarrow \M_{g,n}\] 
is the clutching map attaching the last two marked points to create a loop, then
\begin{equation*}
\gamma^*\mathbb{V}^{\mathfrak{g},\ell}_{g,n}(\mu_1,\ldots,\mu_n) \cong \bigoplus_{\mu\in P_\ell}\mathbb{V}^{\mathfrak{g},\ell}_{g-1,n+2}(\mu_1,\ldots,\mu_n,\mu,\mu^*)
\end{equation*}
and so
\begin{equation}\label{eq:fact2}
\rk\mathbb{V}^{\mathfrak{g},\ell}_{g,n}(\mu_1,\ldots,\mu_n) = \sum_{\mu\in P_\ell}\rk\mathbb{V}^{\mathfrak{g},\ell}_{g-1,n+2}(\mu_1,\ldots,\mu_n,\mu,\mu^*).
\end{equation}

Another ingredient we shall need is the \emph{propagation of vacua} property (also first proven in \cite{TUY89} then expressed in modular terms in \cite[Proposition 2.4]{Fak12}), which says that if a weight $\mu_i$ in $\underline{\mu}$ is zero then 
\begin{equation*}
\mathbb{V}^{\mathfrak{g},\ell}_{g,n}(\underline{\mu}) \cong \pi_i^*\mathbb{V}^{\mathfrak{g},\ell}_{g,n-1}(\mu_1,\ldots,\widehat{\mu_i},\ldots,\mu_n),
\end{equation*}
where $\pi_i$ is the map forgetting the $i^{th}$ marked point.  This says that all points marked by the zero weight can be ignored when computing the rank of a conformal blocks vector bundle.

The final ingredient we shall need is the \emph{fusion rules}, which are computations of the rank when $g=0$ and $n=3$.  The special case of this that we require is the following \cite[Corollary 3.5.2]{Uen97}:  
\begin{equation}\label{eq:2pt}
\rk\mathbb{V}^{\mathfrak{g},\ell}_{0,3}(\mu_1,\mu_2,0) = 
  \begin{cases} 
      1 & \mu_2 = \mu_1^* \\
      0 & \text{otherwise}.
   \end{cases}
\end{equation}

\subsection{Specializing to an exceptional Lie algebra}

Our focus is on the conformal blocks vector bundles for the exceptional Lie algebra $\mathfrak{g} = \mathfrak{g}_2$ and level $\ell = 1$.  Here there is only one nontrivial representation, and it is self-dual, so we can write $P_1 = \{0,\mu\}$ where $\mu^* = \mu$.

The Verlinde formula was proven using factorization, and Mukhopadhyay used factorization and induction to generalize Gregoire's Verlinde evaluation formula \eqref{eq:Gregoire} for $\rk \mathbb{V}^{\mathfrak{g}_2,1}_{g}$ from $n=0$ to all $n\ge 0$ as follows \cite{Muk16}:
\begin{equation}\label{eq:VerG2}
\rk \mathbb{V}^{\mathfrak{g}_2,1}_{g,n}(\mu^{\underline{n}}) = \left(\frac{1+\sqrt{5}}{2}\right)^n\left(\frac{5+\sqrt{5}}{2}\right)^{g-1} + ~  \left(\frac{1-\sqrt{5}}{2}\right)^n\left(\frac{5-\sqrt{5}}{2}\right)^{g-1},
\end{equation}
where \[\mu^{\underline{n}} := (\underbrace{\mu,\ldots,\mu}_{n})\in P_1^{n}.\]  This provides the middle expression in Equation \eqref{eq:main} for all $n \ge 0$ and $g \ge 2$.

\section{Fibonacci formulas from factorization}

\begin{figure}
\begin{center}
(a) \scalebox{0.35}{\includegraphics{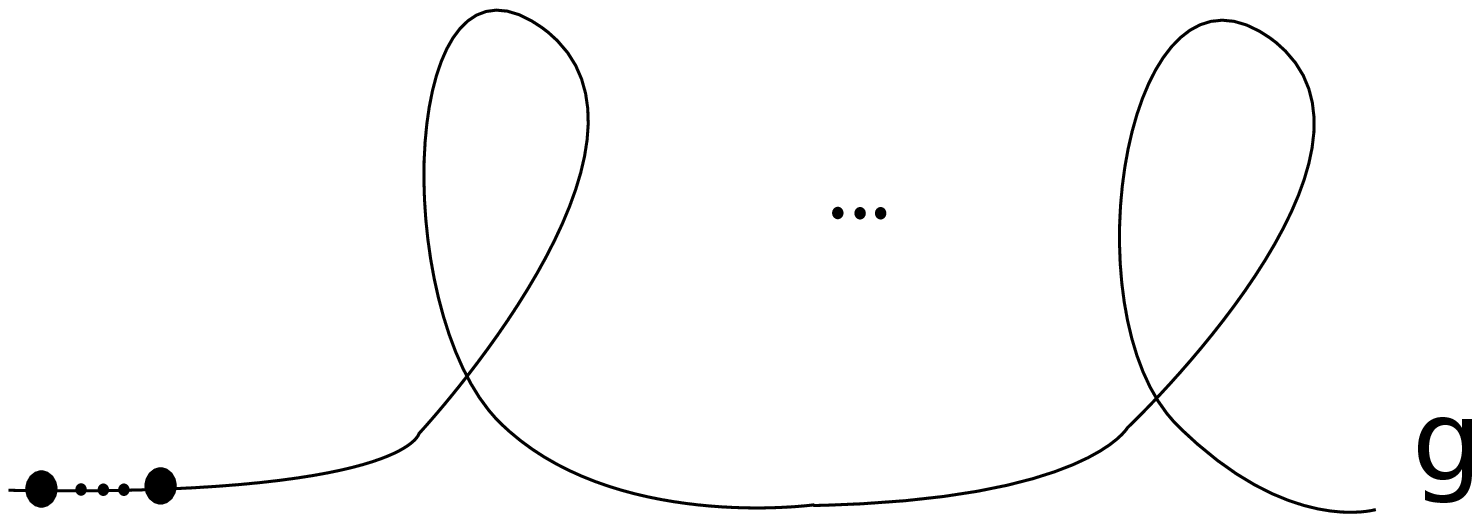}}\\
~\\
(b) \scalebox{0.35}{\includegraphics{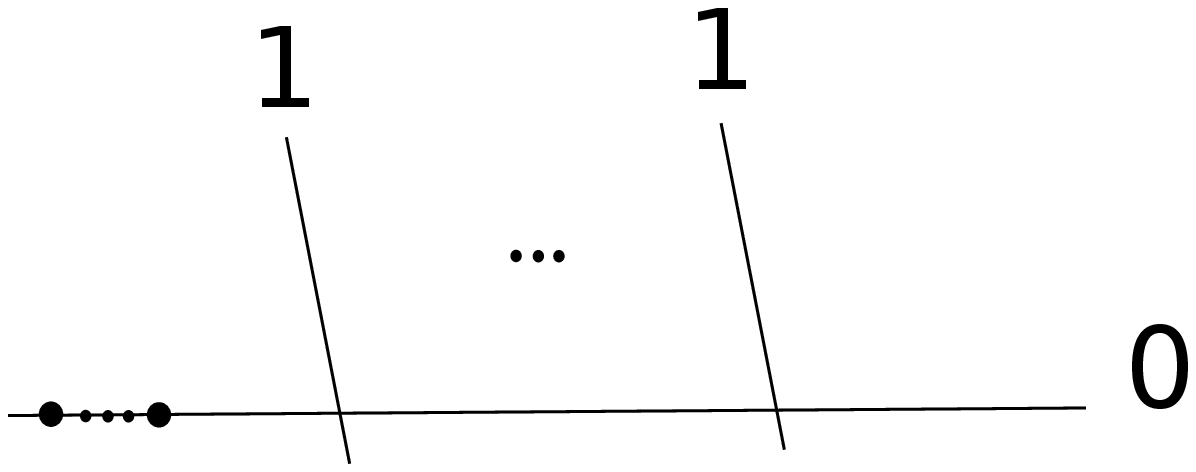}}\\
\end{center}
\caption{(a) An $n$-pointed irreducible genus $g$ curve obtained by clutching $g$ pairs of points on $\mathbb{P}^1$, and (b) a reducible genus $g$ curve obtained by gluing $g$ elliptic tails to an $n$-pointed $\mathbb{P}^1$.}
\label{fig:curve}
\end{figure}

Consider an $n$-pointed genus $g$ curve with $g$ irreducible nodes as in Figure \ref{fig:curve}(a).  This curve represents the generic point of a stratum deep in the stratification of the boundary of $\M_{g,n}$.  A repeated application of factorization (the clutching type described in Equation \eqref{eq:fact2}) yields the following, for any $\underline{\mu}\in P_\ell^n$:
\begin{equation}
\rk\mathbb{V}^{\mathfrak{g},\ell}_{g,n}(\underline{\mu}) = \sum_{(\lambda_1,\ldots,\lambda_g)\in P_\ell^g}\rk\mathbb{V}^{\mathfrak{g},\ell}_{0,n+2g}(\underline{\mu},\lambda_1,\lambda_1^*,\ldots,\lambda_g,\lambda_g^*)
\end{equation}
Specializing to $\mathfrak{g}=\mathfrak{g}_2$ and $\ell=1$, so that $P_1 = \{0,\mu\}$ with $\mu^*=\mu$, we then have
\begin{align} 
\rk\mathbb{V}^{\mathfrak{g}_2,1}_{g,n}(\mu^{\underline{n}}) &= \sum_{(\lambda_1,\ldots,\lambda_g)\in \{0,\mu\}^g}\rk\mathbb{V}^{\mathfrak{g}_2,1}_{0,n+2g}(\mu^{\underline{n}},\lambda_1,\lambda_1,\ldots,\lambda_g,\lambda_g)\\
 &= \sum_{i=0}^g\binom{g}{i}\rk\mathbb{V}^{\mathfrak{g}_2,1}_{0,n+2i}(\mu^{\underline{n+2i}}),
\end{align}
where for the last equality we used propagation of vacua, i.e., weight zero points can be dropped.  Fakhruddin has used factorization, together with the 3-point fusion rule $\rk\mathbb{V}^{\mathfrak{g}_2,1}_{0,3}(\mu,\mu,\mu) = 1$, to show \cite[\S5.2.8]{Fak12} that 
\begin{equation}
\rk\mathbb{V}^{\mathfrak{g}_2,1}_{0,m}(\mu^{\underline{m}}) = F_{m-1},
\end{equation}
the $(m-1)^{th}$ Fibonacci number, and this turns the preceding formula into the left expression in \eqref{eq:main}.

To get the right expression in Equation \eqref{eq:main}, we apply the reducible gluing type of factorization rather than the irreducible clutching applied above.  Specifically, consider the $n$-pointed genus $g$ curve in Figure \ref{fig:curve}(b) where there is a single irreducible ``spine'' component isomorphic to $\mathbb{P}^1$ carrying all $n$ marked points, and attached to this spine are $g$ irreducible elliptic tails.  The factorization formula \eqref{eq:fact1}, for any $\underline{\mu}\in P_\ell^n$, yields
\begin{equation}
\rk\mathbb{V}^{\mathfrak{g},\ell}_{g,n}(\underline{\mu}) = \sum_{(\lambda_1,\ldots,\lambda_g)\in P_\ell^g} \rk\mathbb{V}^{\mathfrak{g},\ell}_{0,n+g}(\underline{\mu},\lambda_1,\ldots,\lambda_g)\prod_{i=1}^g \rk\mathbb{V}^{\mathfrak{g},\ell}_{1,1}(\lambda_i^*).
\end{equation}
Now specializing to $\mathfrak{g}=\mathfrak{g}_2$ and $\ell=1$, so that $P_1 = \{0,\mu\}$ with $\mu^*=\mu$, we have
\begin{align} 
\rk\mathbb{V}^{\mathfrak{g}_2,1}_{g,n}(\mu^{\underline{n}}) &= \sum_{(\lambda_1,\ldots,\lambda_g)\in \{0,\mu\}^g}\rk\mathbb{V}^{\mathfrak{g}_2,1}_{0,n+g}(\mu^{\underline{n}},\lambda_1,\ldots,\lambda_g)\prod_{i=1}^g \rk\mathbb{V}^{\mathfrak{g},\ell}_{1,1}(\lambda_i) \\
 &= \sum_{i=0}^g\binom{g}{i}\rk\mathbb{V}^{\mathfrak{g}_2,1}_{0,n+i}(\mu^{\underline{n+i}})\left(\rk\mathbb{V}^{\mathfrak{g},\ell}_{1,1}(\mu)\right)^i\left(\rk\mathbb{V}^{\mathfrak{g},\ell}_{1,1}(0)\right)^{g-i} \\
  &= \sum_{i=0}^g\binom{g}{i}F_{n+i-1}2^{g-i}.
\end{align}
For the last equality we used that $\rk\mathbb{V}^{\mathfrak{g}_2,1}_{1,1}(\mu) = 1$ and $\rk\mathbb{V}^{\mathfrak{g}_2,1}_{1,1}(0) = 2$, which follow, for instance, from clutching factorization and the fusion rule stated in Equation \eqref{eq:2pt}:
\begin{align*}
\rk\mathbb{V}^{\mathfrak{g}_2,1}_{1,1}(\mu) &= \rk\mathbb{V}^{\mathfrak{g}_2,1}_{0,3}(\mu,0,0) + \rk\mathbb{V}^{\mathfrak{g}_2,1}_{0,3}(\mu,\mu,\mu)\\
 &= 0 + 1,\\
 \rk\mathbb{V}^{\mathfrak{g}_2,1}_{1,1}(0) &= \rk\mathbb{V}^{\mathfrak{g}_2,1}_{0,3}(0,0,0) + \rk\mathbb{V}^{\mathfrak{g}_2,1}_{0,3}(0,\mu,\mu)\\
 &= 1 + 1.
\end{align*}
This completes the proof of Theorem \ref{thm:main}.

\end{document}